\documentclass[12pt]{article}
\usepackage{mathrsfs}
\usepackage{amsfonts}

\usepackage{latexsym}
\usepackage{amssymb}
\usepackage{color}
\usepackage{graphicx}
\usepackage{epstopdf}
\usepackage{amsmath}    
\usepackage{bm}                     
\usepackage{algorithm, algorithmic} 

\newtheorem{Theorem}{Theorem}
\newtheorem{Definition}{Definition}
\newtheorem{Proposition}{Proposition}

\newtheorem{Lemma}{Lemma}

\newtheorem{Example}{Example}

\newcommand{\g}{{\bf g}}
\newcommand{\p}{{\bf g}} 

\newcommand{\A}{{\cal A}}
\newcommand{\B}{{\cal B}}
\newcommand{\C}{{\cal C}}
\newcommand{\D}{{\cal D}}
\newcommand{\HH}{{\cal H}}
\newcommand{\I}{{\cal I}}

\newcommand{\uu}{{\bf u}}
\newcommand{\vv}{{\bf v}}
\newcommand{\x}{{\bf x}}

\newcommand{\y}{{\bf y}}

\renewcommand{\Re}{\mathbb{R}}

\newcommand{\qed}{\nobreak \ifvmode \relax \else
      \ifdim\lastskip<1.5em \hskip-\lastskip
      \hskip1.5em plus0em minus0.5em \fi \nobreak
      \vrule height0.75em width0.5em depth0.25em\fi}

\def \ep{\hbox{ }\hfill$\Box$}

\begin{document}
\title{Computing Eigenvalues of Large Scale Hankel Tensors}

\author{Yannan Chen\footnote{School of Mathematics and Statistics, Zhengzhou University, Zhengzhou, China.
    E-mail: ynchen@zzu.edu.cn (Y. Chen).
    This author's work was supported by the National Natural Science Foundation of China (Grant No. 11401539) and
    the Development Foundation for Excellent Youth Scholars of Zhengzhou University (Grant No. 1421315070).}
    \quad Liqun Qi\footnote{Department of Applied Mathematics, The Hong Kong Polytechnic University,
    Hung Hom, Kowloon, Hong Kong. E-mail: maqilq@polyu.edu.hk (L. Qi).
    This author's work was partially supported by the Hong Kong Research Grant Council
    (Grant No. PolyU 502111, 501212, 501913 and 15302114).}
    \quad Qun Wang\footnote{Department of Applied Mathematics, The Hong Kong Polytechnic University,
    Hung Hom, Kowloon, Hong Kong. Email: wangqun876@gmail.com (Q. Wang).} }

\date{\today} \maketitle

\begin{abstract}

Large scale tensors, including large scale Hankel tensors,
have many applications in science and engineering.
In this paper, we propose an inexact curvilinear search optimization method
to compute Z- and H-eigenvalues of $m$th order $n$ dimensional Hankel tensors, where $n$ is large.
Owing to the fast Fourier transform, the computational cost of each iteration
of the new method is about $\mathcal{O}(mn\log(mn))$.
Using the Cayley transform, we obtain an effective curvilinear search scheme.
Then, we show that every limiting point of iterates generated by the new algorithm is
an eigen-pair of Hankel tensors. Without the assumption of a second-order sufficient condition,
we analyze the linear convergence rate of iterate sequence by the Kurdyka-{\L}ojasiewicz property.
 Finally, numerical experiments for Hankel tensors,
whose dimension may up to one million, are reported to show the efficiency
of the proposed curvilinear search method.

\vspace{3mm}
\noindent {\bf Key words:}\hspace{2mm} Cayley transform, curvilinear search, eigenvalue,
  fast Fourier transform, Hankel tensor, Kurdyka-{\L}ojasiewicz property, large scale tensor.

\vspace{3mm}
\noindent {\bf AMS subject classifications (2010):}\hspace{2mm} 15A18, 15A69, 65F15, 65K05, 90C52.
\vspace{3mm}
\end{abstract}

\newpage
\section{Introduction}

With the coming era of massive data, large scale tensors
have important applications in science and engineering.
How to store and analyze these tensors? This is a pressing and challenging problem.
In the literature, there are two strategies for manipulating large scale tensors.
The first one is to exploit their structures such as sparsity \cite{BK07}.
For example, we consider an online store (e.g. Amazon.com)
where users may review various products \cite{ML13}.
Then, a third order tensor with modes: users, items, and words could be formed naturally
and it is sparse.
The other one is to use distributed and parallel computation \cite{DK14,CiP09}.
This technique could deal with large scale dense tensors, but it depends on a supercomputer.
Recently, researchers applied these two strategies simultaneously for large scale tensors \cite{KPHF12,CV14}.

In this paper, we consider a class of large scale dense tensors with a special Hankel structure.
Hankel tensors appear in many engineering problems such as
signal processing \cite{BDA07,DQW}, automatic control \cite{S14}, and geophysics \cite{OS11,TBM13}.
For instance, in nuclear magnetic resonance spectroscopy \cite{VCDV}, a Hankel matrix was formed
to analyze the time-domain signals, which is important for brain tumour detection.
Papy et al. \cite{PDV05,PDV09} improved this method by using a high order Hankel tensor
to replace the Hankel matrix. Ding et al. \cite{DQW} proposed a fast computational framework
for products of a Hankel tensor and vectors.
On the mathematical properties, Luque and Thibon \cite{LT03} explored the Hankel hyperdeterminants.
Qi \cite{Qi15} and Xu \cite{Xu} studied the spectra of Hankel tensors and
gave some upper bounds and lower bounds for the smallest and the largest eigenvalues.
In \cite{Qi15}, Qi raised a question: Can we construct some efficient algorithms
for the largest and the smallest H- and Z-eigenvalues of a Hankel tensor?

Numerous applications of the eigenvalues of higher order tensors
have been found in science and engineering, such as
automatic control \cite{NQW08}, medical imaging \cite{SS08,QYX13,CDHS13}, quantum information \cite{NQB14},
and spectral graph theory \cite{CD12}.
For example, in magnetic resonance imaging \cite{QYX13}, the principal Z-eigenvalues of an even order tensor
associated to the fiber orientation distribution of a voxel in white matter of human brain
denote volume factions of several nerve fibers in this voxel,
and the corresponding Z-eigenvectors express the orientations of these nerve fibers.
The smallest eigenvalue of tensors reflects the stability of a nonlinear multivariate autonomous system
in automatic control \cite{NQW08}.
For a given even order symmetric tensor, it is positive semidefinite if and only if
its smallest H- or Z-eigenvalue is nonnegative \cite{Qi05}.

The conception of eigenvalues of higher order tensors was defined independently by
Qi \cite{Qi05} and Lim \cite{Lim} in 2005.
Unfortunately, it is an NP-hard problem to compute eigenvalues of a tensor
even though the involved tensor is symmetric \cite{HL13}.
For two and three dimensional symmetric tensors, Qi et al. \cite{QWW09} proposed a direct method
to compute all of its Z-eigenvalues.
It was pointed out in \cite{KM11,KM14} that the polynomial system solver, {\tt NSolve} in \textit{Mathematica},
could be used to compute all of the eigenvalues of lower order and low dimensional tensors.
We note that the mathematical software \textit{Maple} has a similar command {\tt solve}
which is also applicable for the polynomial systems of eigenvalues of tensors.

For general symmetric tensors, Kolda and Mayo \cite{KM11} proposed
a shifted symmetric higher order power method to compute its Z-eigenpairs.
Recently, they \cite{KM14} extended the shifted power method to generalized eigenpairs of tensors
and gave an adaptive shift.
Based on the nonlinear optimization model with a compact unit spherical constraint,
the power methods \cite{DDV00} project the gradient of the objective at the current iterate
onto the unit sphere at each iteration.
Its computation is very simple but may not converge \cite{KR02}.
Kolda and Mayo \cite{KM11,KM14} introduced a shift to force the objective to be (locally) concave/convex.
Then the power method produces increasing/decreasing steps for computing maximal/minimal eigenvalues.
The sequence of objectives converges to eigenvalues since the feasible region is compact.
The convergence of the sequence of iterates to eigenvectors is established under the assumption that
the tensor has finitely many real eigenvectors.
The linear convergence rate is estimated by a fixed-point analysis.

Inspired by the power method, various optimization methods have been established.
Han \cite{Han13} proposed an unconstrained optimization model, which is indeed a quadratic penalty function
of the constrained optimization for generalized eigenvalues of symmetric tensors.
Hao et al. \cite{HCD} employed a subspace projection method for Z-eigenvalues of symmetric tensors.
Restricted by a unit spherical constraint, this method minimizes the objective in
a big circle of $n$ dimensional unit sphere at each iteration.
Since the objective is a homogeneous polynomial, the minimization of the subproblem has a closed-form solution.
Additionally, Hao et al. \cite{HCD15a} gave a trust region method to calculate Z-eigenvalues
of symmetric tensors. The sequence of iterates generated by this method converges to
a second order critical point and enjoys a locally quadratic convergence rate.

Since nonlinear optimization methods may produce a local minimizer,
some convex optimization models have been studied.
Hu et al. \cite{HHQ} addressed a sequential semi-definite programming method to compute the extremal
Z-eigenvalues of tensors. A sophisticated Jacobian semi-definite relaxation method was explored
by Cui et al. \cite{CDN}. A remarkable feature of this method is the ability to compute
all of the real eigenvalues of symmetric tensors.
Recently, Chen et al. \cite{CHZ} proposed homotopy continuation methods to
compute all of the complex eigenvalues of tensors.
When the order or the dimension of a tensor grows larger,
the CPU times of these methods become longer and longer.

In some applications \cite{VCDV,OS11}, the scales of Hankel tensors can be quite huge.
This highly restricted the applications of the above mentioned methods in this case.
How to compute the smallest and the largest eigenvalues of a Hankel tensor?
Can we propose a method to compute the smallest and the largest eigenvalues of a
relatively large Hankel tensor, say $1,000,000$ dimension?
This is one of the motivations of this paper.

Owing to the multi-linearity of tensors, we model the problem of eigenvalues of Hankel tensors
as a nonlinear optimization problem with a unit spherical constraint.
Our algorithm is an inexact steepest descent method on the unit sphere.
To preserve iterates on the unit sphere, we employ the Cayley transform
to generate an orthogonal matrix such that the new iterate is this orthogonal matrix times the current iterate.
By the Sherman-Morrison-Woodbury formula, the product of the orthogonal matrix and a vector
has a closed-form solution. So the subproblem is straightforward.
A curvilinear search is employed to guarantee the convergence.
Then, we prove that every accumulation point of the sequence of iterates is an eigenvector of
the involved Hankel tensor, and its objective is the corresponding eigenvalue.
Furthermore, using the Kurdyka-{\L}ojasiewicz property of the eigen-problem of tensors,
we prove that the sequence of iterates converges without
an assumption of second order sufficient condition.
Under mild conditions, we show that the sequence of iterates has
a linear or a sublinear convergence rate.
Numerical experiments show that this strategy is successful.


The outline of this paper is drawn as follows.
We introduce a fast computational framework for products of a well-structured Hankel tensor
and vectors in Section 2. The computational cost is cheap.
In Section 3, we show the techniques of using the Cayley transform
to construct an effective curvilinear search algorithm.
The convergence of objective and iterates are analyzed in Section 4.
The Kurdyka-{\L}ojasiewicz property is applied to analyze an inexact line search method.
Numerical experiments in Section 5 address that the new method is efficient and promising.
Finally, we conclude the paper with Section 6.

\section{Hankel tensors}

Suppose $\A$ is an $m$th order $n$ dimensional real symmetric tensor
\begin{equation*}
    \A = (a_{i_1,i_2,\ldots,i_m}), \qquad \text{ for }i_j=1,\ldots,n, j=1,\ldots,m,
\end{equation*}
where all of the entries are real and invariant under any index permutation.
Two products of the tensor $\A$ and a column vector $\x\in\Re^n$ used in this paper are defined as follows.
\begin{itemize}
  \item $\A\x^m$ is a scalar
    \begin{equation*}
        \A\x^m = \sum_{i_1,\ldots,i_m=1}^n a_{i_1,\ldots,i_m}x_{i_1}\cdots x_{i_m}.
    \end{equation*}
  \item $\A\x^{m-1}$ is a column vector
    \begin{equation*}
        \left(\A\x^{m-1}\right)_i = \sum_{i_2,\ldots,i_m=1}^n a_{i,i_2,\ldots,i_m}x_{i_2}\cdots x_{i_m},
          \qquad \text{for }i=1,\ldots,n.
    \end{equation*}
\end{itemize}
When the tensor $\A$ is dense,
the computations of produces $\A\x^m$ and $\A\x^{m-1}$ require $\mathcal{O}(n^m)$ operations,
since the tensor $\A$ has $n^m$ entries and we must visit all of them in the process of calculation.
When the tensor is symmetric, the computational cost for these products is about
$\mathcal{O}(n^m/m!)$ \cite{SLVK14}. Obviously, they are expensive.
In this section, we will study a special tensor, the Hankel tensor, whose elements are completely
determined by a generating vector. So there exists a fast algorithm to compute products of
a Hankel tensor and vectors. Let us give the definitions of two structured tensors.

\begin{Definition}
  An $m$th order $n$ dimensional tensor $\HH$ is called a Hankel tensor if its entries satisfy
  \begin{equation*}
    h_{i_1,i_2,\ldots,i_m} = v_{i_1+i_2+\cdots+i_m-m}, \qquad \text{for }i_j=1,\ldots,n, j=1,\ldots,m.
  \end{equation*}
  The vector $\vv=(v_0,v_1,\ldots,v_{m(n-1)})^{\top}$ with length $\ell\equiv m(n-1)+1$ is called the generating vector of the Hankel tensor $\HH$.

  An $m$th order $\ell$ dimensional tensor $\C$ is called an anti-circulant tensor if its entries satisfy
  \begin{equation*}
    c_{i_1,i_2,\ldots,i_m} = v_{(i_1+i_2+\cdots+i_m-m~\mathrm{mod}~\ell)},
      \qquad \text{for }i_j=1,\ldots,\ell, j=1,\ldots,m.
  \end{equation*}
\end{Definition}

It is easy to see that $\HH$ is a sub-tensor of $\C$.
Since for the same generating vector $\vv$ we have
\begin{equation*}
    c_{i_1,i_2,\ldots,i_m} = h_{i_1,i_2,\ldots,i_m},
      \qquad \text{for }i_j=1,\ldots,n, j=1,\ldots,m.
\end{equation*}
For example, a third order two dimensional Hankel tensor
with a generating vector $\vv=(v_0,v_1,v_2,v_3)^{\top}$ is
\begin{equation*}
    \HH = \left[
            \begin{array}{cc|cc}
              v_0 & v_1 & v_1 & v_2 \\
              v_1 & v_2 & v_2 & v_3 \\
            \end{array}
          \right].
\end{equation*}
It is a sub-tensor of an anti-circulant tensor with the same order and a larger dimension
\begin{equation*}
    \C = \left[
            \begin{array}{cccc|cccc|cccc|cccc}
              v_0 & v_1 & v_2 & v_3 & v_1 & v_2 & v_3 & v_0 & v_2 & v_3 & v_0 & v_1 & v_3 & v_0 & v_1 & v_2 \\
              v_1 & v_2 & v_3 & v_0 & v_2 & v_3 & v_0 & v_1 & v_3 & v_0 & v_1 & v_2 & v_0 & v_1 & v_2 & v_3 \\
              v_2 & v_3 & v_0 & v_1 & v_3 & v_0 & v_1 & v_2 & v_0 & v_1 & v_2 & v_3 & v_1 & v_2 & v_3 & v_0 \\
              v_3 & v_0 & v_1 & v_2 & v_0 & v_1 & v_2 & v_3 & v_1 & v_2 & v_3 & v_0 & v_2 & v_3 & v_0 & v_1 \\
            \end{array}
          \right].
\end{equation*}

As discovered in \cite[Theorem 3.1]{DQW}, the $m$th order $\ell$ dimensional anti-circulant tensor $\C$
could be diagonalized by the $\ell$-by-$\ell$ Fourier matrix $F_{\ell}$, i.e., $\C = \D F_{\ell}^m$,
where $\D$ is a diagonal tensor whose diagonal entries are $\mathrm{diag}(\D) = F_{\ell}^{-1}\vv$.
It is well-known that the computations involving the Fourier matrix and its inverse times a vector
are indeed the fast (inverse) Fourier transform {\tt fft} and {\tt ifft}, respectively.
The computational cost is about $\mathcal{O}(\ell\log\ell)$ multiplications,
which is significantly smaller than $\mathcal{O}(\ell^2)$ for a dense matrix times a vector
when the dimension $\ell$ is large.


Now, we are ready to show how to compute the products introduced in the beginning of this section,
when the involved tensor has a Hankel structure.
For any $\x\in\Re^n$, we define another vector $\y\in\Re^{\ell}$ such that
\begin{equation*}
    \y \equiv \left[
                \begin{array}{c}
                  \x \\
                  {\bf 0}_{\ell-n} \\
                \end{array}
              \right],
\end{equation*}
where $\ell=m(n-1)+1$ and ${\bf 0}_{\ell-n}$ is a zero vector with length $\ell-n$.
Then, we have
\begin{equation*}
  \HH\x^m = \C\y^m = \D (F_{\ell}\y)^m
    = \mathrm{\text{\tt ifft}}(\vv)^{\top} \left(\mathrm{\text{\tt fft}}(\y)^{\circ m}\right).
\end{equation*}
To obtain $\HH\x^{m-1}$, we first compute
\begin{equation*}
    \C\y^{m-1} = F_{\ell} \left(\D (F_{\ell}\y)^{m-1}\right)
    = \mathrm{\text{\tt fft}}\left(
        \mathrm{\text{\tt ifft}}(\vv) \circ \left(\mathrm{\text{\tt fft}}(\y)^{\circ (m-1)}\right)
      \right).
\end{equation*}
Then, the entries of vector $\HH\x^{m-1}$ is the leading $n$ entries of $\C\y^{m-1}$.
Here, $\circ$ denotes the Hadamard product such that $(A\circ B)_{i,j}=A_{i,j}B_{i,j}$.
Three matrices $A$, $B$ and $A\circ B$ have the same size.
Furthermore, we define $A^{\circ k}=A\circ \cdots\circ A$ as the Hadamard product of
$k$ copies of $A$.

Since the computations of $\HH\x^m$ and $\HH\x^{m-1}$ require $2$ and $3$ {\tt fft/ifft}s,
the computational cost is about $\mathcal{O}(mn\log(mn))$ and obviously cheap.
Another advantage of this approach is that we do not need to store and deal with
the tremendous Hankel tensor explicitly. It is sufficient to keep and work with
the compact generating vector of that Hankel tensor.


\section{A curvilinear search algorithm}

\newcommand{\SPHERE}{\mathbb{S}_{n-1}}

We consider the generalized eigenvalue \cite{CPZ09,DW15} of an $m$th order $n$ dimensional Hankel tensor $\HH$
\begin{equation*}
    \HH\x^{m-1} = \lambda\B\x^{m-1},
\end{equation*}
where $m$ is even, $\B$ is an $m$th order $n$ dimensional symmetric tensor
and it is positive definite.
If there is a scalar $\lambda$ and a real vector $\x$ satisfying this system,
we call $\lambda$ a generalized eigenvalue and $\x$ its associated generalized eigenvector.
Particularly, we find the following definitions from the literature,
where the computation on the tensor $\B$ is straightforward.
\begin{itemize}
  \item Qi \cite{Qi05} called a real scalar $\lambda$ a Z-eigenvalue of a tensor $\HH$ and
    a real vector $\x$ its associated Z-eigenvector if they satisfy
    \begin{equation*}
        \HH\x^{m-1}=\lambda\x ~~~~\text{and}~~~~\x^{\top}\x=1.
    \end{equation*}
    This definition means that the tensor $\B$ is an identity tensor ${\cal E}$
    such that ${\cal E}\x^{m-1}=\|\x\|^{m-2}\x$.
  \item If $\B=\I$, where
    \begin{equation*}
    (\I)_{i_1,\ldots,i_m}=\left\{\begin{aligned}
      1  &~~~~\text{if }i_1=\cdots=i_m, \\
      0  &~~~~\text{ otherwise },
    \end{aligned}\right.
    \end{equation*}
    the real scalar $\lambda$ is called an H-eigenvalue
    and the real vector $\x$ is its associated H-eigenvector \cite{Qi05}.
    Obviously, we have $(\I\x^{m-1})_i=x_i^{m-1}$ for $i=1,\ldots,n$.
\end{itemize}

To compute a generalized eigenvalue and its associated eigenvector,
we consider the following optimization model
with a spherical constraint
\begin{equation}\label{eq3-02}
    \min~~f(\x) \equiv \frac{\HH\x^m}{\B\x^m} \qquad
    \mathrm{s.t.}~~ \|\x\|=1,
\end{equation}
where $\|\cdot\|$ denotes the Euclidean norm or its induced matrix norm.
The denominator of the objective is positive since the tensor $\B$ is positive definite.
By some calculations, we get its gradient and Hessian, which are formally presented in the following lemma.

\begin{Lemma}\label{lm3-1}
  Suppose that the objective is defined as in (\ref{eq3-02}). Then, its gradient is
  \begin{equation}\label{eq3-03}
    \g(\x)=\frac{m}{\B\x^m}\left(\HH\x^{m-1}-\frac{\HH\x^m}{\B\x^m}\B\x^{m-1}\right).
  \end{equation}
  And its Hessian is
  \begin{eqnarray}\label{eq3-04}
    H(\x)&=&\frac{m(m-1)\HH\x^{m-2}}{\B\x^m}
      -\frac{m(m-1)\HH\x^m\B\x^{m-2}+m^2(\HH\x^{m-1}\circledcirc\B\x^{m-1})}{(\B\x^m)^2} \nonumber\\
      &&{} +\frac{m^2\HH\x^m(\B\x^{m-1}\circledcirc\B\x^{m-1})}{(\B\x^m)^3},
  \end{eqnarray}
  where $\x\circledcirc\y \equiv \x\y^{\top}+\y\x^{\top}$.
\end{Lemma}


Let $\SPHERE\equiv\{\x\in\Re^n~|~\x^{\top}\x=1\}$ be the spherical feasible region.
Suppose the current iterate is $\x\in\SPHERE$ and the gradient at $\x$ is $\g(\x)$.
Because
\begin{equation}\label{eq3-06}
    \x^{\top}\g(\x) = \frac{m}{\B\x^m}\left(\x^{\top}\HH\x^{m-1}
            -\frac{\HH\x^m}{\B\x^m}\x^{\top}\B\x^{m-1}\right)  = 0,
\end{equation}
the gradient $\g(\x)$ of $\x\in\SPHERE$ is located in the tangent plane of $\SPHERE$ at $\x$.

%

\begin{Lemma}\label{Lm3-2}
  Suppose $\|\g(\x)\|=\epsilon$, where $\x\in\SPHERE$ and $\epsilon$ is a small number.
  Denote $\lambda=\frac{\HH\x^m}{\B\x^m}$.
  Then, we have
  \begin{equation*}
    \| \HH\x^{m-1}-\lambda\B\x^{m-1} \| = \mathcal{O}(\epsilon).
  \end{equation*}
  Moreover, if the gradient $\g(\x)$ at $\x$ vanishes, then
  $\lambda=f(\x)$ is a generalized eigenvalue
  and $\x$ is its associated generalized eigenvector.
\end{Lemma}
\noindent
{\bf Proof}
  Recalling the definition of gradient (\ref{eq3-03}), we have
  \begin{equation*}
    \| \HH\x^{m-1}-\lambda\B\x^{m-1} \| = \frac{\B\x^m}{m}\epsilon.
  \end{equation*}
  Since the tensor $\B$ is positive definite and the vector $\x$ belongs to a compact set $\SPHERE$,
  $\B\x^m$ has a finite upper bound. Thus, the first assertion is valid.

  If $\epsilon =0$, we immediately know that $\lambda=f(\x)$ is a generalized eigenvalue
  and $\x$ is its associated generalized eigenvector.
\ep

Next, we construct the curvilinear search path using the Cayley transform \cite{GV13}.
Cayley transform is an effective method which could preserve the orthogonal constraints.
It has various applications in the inverse eigenvalue problem \cite{FNO87}, $p$-harmonic flow \cite{GWY09},
and matrix optimization \cite{WY}.

Suppose the current iterate is $\x_k\in\SPHERE$ and the next iterate is $\x_{k+1}$.
To preserve the spherical constraint $\x_{k+1}^{\top}\x_{k+1}=\x_k^{\top}\x_k=1$, we choose
the next iterate $\x_{k+1}$ such that
\begin{equation}\label{eq3-07}
    \x_{k+1} = Q\x_k,
\end{equation}
where $Q\in\Re^{n \times n}$ is an orthogonal matrix,
whose eigenvalues do not contain $-1$.
Using the Cayley transform, the matrix
\begin{equation}\label{eq3-08}
    Q=(I+W)^{-1}(I-W)
\end{equation}
is orthogonal if and only if the matrix $W\in\Re^{n\times n}$ is
skew-symmetric.\footnote{See ``http://en.wikipedia.org/wiki/Cayley\_transform".}
Now, our task is to select a suitable skew-symmetric matrix $W$
such that $\p(\x_k)^{\top}(\x_{k+1}-\x_k) < 0$.
For simplicity, we take the matrix $W$ as
\begin{equation}\label{eq3-09}
    W={\bf a}{\bf b}^{\top} - {\bf b}{\bf a}^{\top},
\end{equation}
where ${\bf a},{\bf b}\in\Re^n$ are two undetermined vectors.
From (\ref{eq3-07}) and (\ref{eq3-08}), we have
\begin{equation*}
    \x_{k+1}-\x_k = -W(\x_k+\x_{k+1}).
\end{equation*}
Then, by (\ref{eq3-09}), it yields that
\begin{equation*}
  \p(\x_k)^{\top}(\x_{k+1}-\x_k) = -[(\p(\x_k)^{\top}{\bf a}){\bf b}^{\top}
    - (\p(\x_k)^{\top}{\bf b}){\bf a}^{\top}](\x_k+\x_{k+1}).
\end{equation*}
For convenience, we choose
\begin{equation}\label{eq3-10}
    {\bf a}=\x_k~~~~\text{and}~~~~{\bf b}=-\alpha\p(\x_k).
\end{equation}
Here, $\alpha$ is a positive parameter, which serves as a step size,
so that we have some freedom to choose the next iterate.
According to this selection and (\ref{eq3-06}), we obtain
\begin{eqnarray*}
    \p(\x_k)^{\top}(\x_{k+1}-\x_k) &=& -\alpha\|\p(\x_k)\|^2 \x_k^{\top}(\x_k+\x_{k+1}) \\
      &=& -\alpha\|\p(\x_k)\|^2 (1+\x_k^{\top}Q\x_k).
\end{eqnarray*}
Since $-1$ is not an eigenvalue of the orthogonal matrix $Q$,
we have $1+\x_k^{\top}Q\x_k > 0$ for $\x_k^{\top}\x_k=1$.
Therefore, the conclusion $\p(\x_k)^{\top}(\x_{k+1}-\x_k)<0$ holds for any positive step size $\alpha$.

We summarize the iterative process in the following Theorem.

\begin{Theorem}
  Suppose that the new iterate $\x_{k+1}$ is generated by (\ref{eq3-07}), (\ref{eq3-08}), (\ref{eq3-09}),
  and (\ref{eq3-10}). Then, the following assertions hold.
  \begin{itemize}
    \item The iterative scheme is
      \begin{equation}\label{eq3-11}
        \x_{k+1}(\alpha) = \frac{1-\alpha^2\|\p(\x_k)\|^2}{1+\alpha^2\|\p(\x_k)\|^2}\x_k
          -\frac{2\alpha}{1+\alpha^2\|\p(\x_k)\|^2}\p(\x_k).
      \end{equation}
    \item The progress made by $\x_{k+1}$ is
      \begin{equation}\label{eq3-12}
        \p(\x_k)^{\top}(\x_{k+1}(\alpha)-\x_k) = -\frac{2\alpha\|\p(\x_k)\|^2}{1+\alpha^2\|\p(\x_k)\|^2}.
      \end{equation}
  \end{itemize}
\end{Theorem}
\noindent
{\bf Proof}
  From the equality (\ref{eq3-06}) and the Sherman-Morrison-Woodbury formula, we have
  \begin{eqnarray*}
    \x_{k+1}(\alpha) &=& (I-\alpha\x_k\p(\x_k)^{\top}+\alpha\p(\x_k)\x_k^{\top})^{-1}
        (I+\alpha\x_k\p(\x_k)^{\top}-\alpha\p(\x_k)\x_k^{\top})\x_k \\
      &=& (I+\alpha\p(\x_k)\x_k^{\top}-\alpha\x_k\p(\x_k)^{\top})^{-1}(\x_k-\alpha\p(\x_k)) \\
      &=& \Bigg(I-\left[
                     \begin{array}{cc}
                       \alpha\p(\x_k) & -\x_k \\
                     \end{array}
                   \right]\left(\left[
                                  \begin{array}{cc}
                                    1 & 0 \\
                                    0 & 1 \\
                                  \end{array}
                                \right]+\left[
                                          \begin{array}{c}
                                            \x_k^{\top} \\
                                            \alpha\p(\x_k)^{\top} \\
                                          \end{array}
                                        \right]I\left[
                                                  \begin{array}{cc}
                                                    \alpha\p(\x_k) & -\x_k \\
                                                  \end{array}
                                                \right]
                   \right)^{-1}\cdot \\
    && {}~~~~~~~~~~~~~~~~~~~~ \left[
                                 \begin{array}{c}
                                   \x_k^{\top} \\
                                   \alpha\p(\x_k)^{\top} \\
                                 \end{array}
                               \right]
      \Bigg)(\x_k-\alpha\p(\x_k)) \\
    &=& \x_k-\alpha\p(\x_k) - \left[
                     \begin{array}{cc}
                       \alpha\p(\x_k) & -\x_k \\
                     \end{array}
                   \right]\left[
                            \begin{array}{cc}
                              1 & -1 \\
                              \alpha^2\|\p(\x_k)\|^2 & 1  \\
                            \end{array}
                          \right]^{-1}\left[
                                        \begin{array}{c}
                                          1 \\
                                          -\alpha^2\|\p(\x_k)\|^2 \\
                                        \end{array}
                                      \right] \\
    &=& \frac{1-\alpha^2\|\p(\x_k)\|^2}{1+\alpha^2\|\p(\x_k)\|^2}\x_k
          -\frac{2\alpha}{1+\alpha^2\|\p(\x_k)\|^2}\p(\x_k).
  \end{eqnarray*}
  The proof of (\ref{eq3-12}) is straightforward.
\ep

Whereafter, we devote to choose a suitable step size $\alpha$ by an inexact curvilinear search.
At the beginning, we give a useful theorem.

\begin{Theorem}\label{th3-5}
  Suppose that the new iterate $\x_{k+1}(\alpha)$ is generated by (\ref{eq3-11}).
  Then, we have
  \begin{equation*}
    \left.\frac{\mathrm{d}f(\x_{k+1}(\alpha))}{\mathrm{d}\alpha}\right|_{\alpha=0} = -2\|\p(\x_k)\|^2.
  \end{equation*}
\end{Theorem}
\noindent
{\bf Proof}
  By some calculations, we get
  \begin{equation*}
    \x_{k+1}'(\alpha) = \frac{-2}{1+\alpha^2\|\p(\x_k)\|^2}\p(\x_k)
      + \frac{-4\alpha\|\p(\x_k)\|^2}{(1+\alpha^2\|\p(\x_k)\|^2)^2}(\x_k-\alpha\p(\x_k)).
  \end{equation*}
  Hence, $\x_{k+1}'(0)=-2\p(\x_k)$. Furthermore, $\x_{k+1}(0)=\x_k$.
  Therefore, we obtain
  \begin{equation*}
    \left.\frac{\mathrm{d}f(\x_{k+1}(\alpha))}{\mathrm{d}\alpha}\right|_{\alpha=0}
      = \g(\x_{k+1}(0))^{\top}\x_{k+1}'(\alpha) = \g(\x_k)^{\top}(-2\p(\x_k))
      = -2\|\p(\x_k)\|^2.
  \end{equation*}
  The proof is completed.
\ep

\begin{algorithm}[tb!]
\caption{A curvilinear search algorithm (ACSA).}\label{alg}
\begin{algorithmic}[1]
  \STATE Give the generating vector $\vv$ of a Hankel tensor $\HH$, the symmetric tensor $\B$,
    an initial unit iterate $\x_1$, parameters $\eta \in (0,\frac{1}{2}]$, $\beta\in(0,1)$,
    $\bar{\alpha}_1=1 \leq \alpha_{\max}$,
    and $k \gets 1$.

  \WHILE{the sequence of iterates does not converge}

    \STATE Compute $\HH\x_k^m$ and $\HH\x_k^{m-1}$ by the fast computational framework introduces in Section 2.
    \STATE Calculate $\B\x_k^m$, $\B\x_k^{m-1}$, $\lambda_k=f(\x_k)=\frac{\HH\x_k^m}{\B\x_k^m}$
       and $\p(\x_k)$ by (\ref{eq3-03}).

    \STATE Choose the smallest nonnegative integer $\ell$ and determine $\alpha_k=\beta^{\ell}\bar{\alpha}_k$
      such that
      \begin{equation}\label{suf-dec}
        f(\x_{k+1}(\alpha_k)) \leq f(\x_k) - \eta\alpha_k\|\p(\x_k)\|^2,
      \end{equation}
      where $\x_{k+1}(\alpha)$ is calculated by (\ref{eq3-11}).


    \STATE Update the iterate $\x_{k+1}=\x_{k+1}(\alpha_k)$.
    \STATE Choose an initial step size $\bar{\alpha}_{k+1}\in(0,\alpha_{\max}]$ for the next iteration.
    \STATE $k \gets k+1.$
  \ENDWHILE
\end{algorithmic}
\end{algorithm}

According to Theorem \ref{th3-5}, for any constant $\eta\in(0,2)$,
there exists a positive scalar $\tilde{\alpha}$ such that for all $\alpha\in(0,\tilde{\alpha}]$,
\begin{equation*}
    f(\x_{k+1}(\alpha)) - f(\x_k) \leq -\eta\alpha\|\p(\x_k)\|^2.
\end{equation*}
Hence, the curvilinear search process is well-defined.

Now, we present a curvilinear search algorithm (ACSA) formally in Algorithm \ref{alg}
for the smallest generalized eigenvalue and its associated eigenvector of a Hankel tensor.
If our aim is to compute the largest generalized eigenvalue and its associated eigenvector
of a Hankel tensor, we only need to change respectively (\ref{eq3-11}) and (\ref{suf-dec})
used in Steps 5 and 6 of the ACSA algorithm to
\begin{equation*}
    \x_{k+1}(\alpha) = \frac{1-\alpha^2\|\p(\x_k)\|^2}{1+\alpha^2\|\p(\x_k)\|^2}\x_k
          +\frac{2\alpha}{1+\alpha^2\|\p(\x_k)\|^2}\p(\x_k),
\end{equation*}
and
\begin{equation*}
    f(\x_{k+1}(\alpha_k)) \geq f(\x_k) + \eta\alpha_k\|\p(\x_k)\|^2.
\end{equation*}

When the Z-eigenvalue of a Hankel tensor is considered,
we have ${\cal E}\x^m=\|\x\|^m=1$ and the objective $f(\x)$ is a polynomial.
Then, we could compute the global minimizer of the step size $\alpha_k$ (the exact line search)
in each iteration as \cite{HCD}.
However, we use a cheaper inexact line search here.
The initial step size of the next iteration follows Dai's strategy \cite{Dai14}
\begin{equation*}
    \bar{\alpha}_{k+1}=\frac{\|\Delta\x_k\|}{\|\Delta\p_k\|},
\end{equation*}
which is the geometric mean of Barzilai-Borwein step sizes \cite{BB88}.


\section{Convergence analysis}

Since the optimization model (\ref{eq3-02}) has a nice algebraic nature,
we will use the Kurdyka-{\L}ojasiewicz property \cite{L63,BDL07} to analyze
the convergence of the proposed ACSA algorithm. Before we start, we give some
basic convergence results.

\subsection{Basic convergence results}

If the ACSA algorithm terminates finitely,
there exists a positive integer $k$ such that $\p(\x_k)=0$.
According to Lemma \ref{Lm3-2}, $f(\x_k)$ is a generalized eigenvalue and
$\x_k$ is its associated generalized eigenvector.

Next, we assume that ACSA generates an infinitely sequence of iterates.

\begin{Lemma}
  Suppose that the even order symmetric tensor $\B$ is positive definite.
  Then, all the functions, gradients, and Hessians of the objective (\ref{eq3-02})
  at feasible points are bounded.
  That is to say, there is a positive constant $M$ such that for all $\x\in\SPHERE$
  \begin{equation}\label{eq4-01}
    |f(\x)|\leq M,\quad \|\p(\x)\|\leq M, \quad\text{and}\quad \|H(\x)\|\leq M.
  \end{equation}
\end{Lemma}
\noindent
{\bf Proof}
  Since the spherical feasible region $\SPHERE$ is compact,
  the denominator $\B\x^m$ of the objective is positive and bounds away from zero.
  Recalling Lemma \ref{lm3-1}, we get this theorem immediately.
\ep

\begin{Theorem}
  Suppose that the infinite sequence $\{\lambda_k\}$ is generated by ACSA.
  Then, the sequence $\{\lambda_k\}$ is monotonously decreasing.
  And there exists a $\lambda_*$ such that
  \begin{equation*}
    \lim_{k\to\infty} \lambda_k=\lambda_*.
  \end{equation*}
\end{Theorem}
\noindent
{\bf Proof}
  Since $\lambda_k=f(\x_k)$ which is bounded and monotonously decreasing,
  the infinite sequence $\{\lambda_k\}$ must converge to a unique $\lambda_*$.
\ep

This theorem means that the sequence of generalized eigenvalues converges.
To show the convergence of iterates, we first prove that the step sizes bound
away from zero.

\begin{Lemma}
  Suppose that the step size $\alpha_k$ is generated by ACSA.
  Then, for all iterations $k$, we get
  \begin{equation}\label{eq4-06}
    \alpha_k \geq \frac{(2-\eta)\beta}{5M}\equiv\alpha_{\min} > 0.
  \end{equation}
\end{Lemma}
\noindent
{\bf Proof}
  Let $\underline{\alpha}\equiv\frac{(2-\eta)}{5M}$.
  According to the curvilinear search process of ACSA, it is sufficient to prove that
  the inequality (\ref{suf-dec}) holds if $\alpha_k\in(0,\underline{\alpha}]$.

  From the iterative formula (\ref{eq3-11}) and the equality (\ref{eq3-06}), we get
  \begin{eqnarray*}
    \|\x_{k+1}(\alpha)-\x_k\|^2 &=& \left\|\frac{-2\alpha^2\|\p(\x_k)\|^2}{1+\alpha^2\|\p(\x_k)\|^2}\x_k
          -\frac{2\alpha}{1+\alpha^2\|\p(\x_k)\|^2}\p(\x_k)\right\|^2 \\
      &=& \frac{4\alpha^4\|\p(\x_k)\|^4\|\x_k\|^2 + 4\alpha^2\|\p(\x_k)\|^2}{(1+\alpha^2\|\p(\x_k)\|^2)^2} \\
      &=& \frac{4\alpha^2\|\p(\x_k)\|^2}{1+\alpha^2\|\p(\x_k)\|^2}.
  \end{eqnarray*}
  Hence,
  \begin{equation}\label{eq4-07}
    \|\x_{k+1}(\alpha)-\x_k\| = \frac{2\alpha\|\p(\x_k)\|}{\sqrt{1+\alpha^2\|\p(\x_k)\|^2}}.
  \end{equation}

  From the mean value theorem, (\ref{eq3-11}), (\ref{eq3-06}), and (\ref{eq4-07}), we have
  \begin{eqnarray*}
    \lefteqn{ f(\x_{k+1}(\alpha))-f(\x_k)
      \leq \g(\x_k)^{\top}(\x_{k+1}(\alpha)-\x_k)+\frac{1}{2}M\|\x_{k+1}(\alpha)-\x_k\|^2 }\\
      &=& \frac{1}{1+\alpha^2\|\p(\x_k)\|^2}\left(-2\alpha^2\|\p(\x_k)\|^2\g(\x_k)^{\top}\x_k
            -2\alpha\|\p(\x_k)\|^2+\frac{M}{2}4\alpha^2\|\p(\x_k)\|^2\right) \\
      &\leq& \frac{\alpha\|\p(\x_k)\|^2}{1+\alpha^2\|\p(\x_k)\|^2}\left(
            4\alpha M-2 \right).
  \end{eqnarray*}
  It is easy to show that for all $\alpha\in(0,\underline{\alpha}]$
  \begin{equation*}
    4\alpha M-2 \leq -\eta(1+\alpha^2M^2).
  \end{equation*}
  Therefore, we have
  \begin{equation*}
    f(\x_{k+1}(\alpha))-f(\x_k) \leq \frac{-\eta(1+\alpha^2M^2)}{1+\alpha^2\|\p(\x_k)\|^2}\alpha\|\p(\x_k)\|^2
      \leq -\eta\alpha\|\p(\x_k)\|^2.
  \end{equation*}
  The proof is completed.
\ep

\begin{Theorem}\label{Th4-04}
  Suppose that the infinite sequence $\{\x_k\}$ is generated by ACSA.
  Then, the sequence $\{\x_k\}$ has an accumulation point at least. And we have
  \begin{equation}\label{eq4-05}
    \lim_{k\to\infty} \|\p(\x_k)\|=0.
  \end{equation}
  That is to say, every accumulation point of $\{\x_k\}$ is a generalized eigenvector
  whose associated generalized eigenvalue is $\lambda_*$.
\end{Theorem}
\noindent
{\bf Proof}
  Since the sequence of objectives $\{f(\x_k)\}$ is monotonously decreasing and bounded,
  by (\ref{suf-dec}) and (\ref{eq4-06}), we have
  \begin{equation*}
    2M \geq f(\x_1)-\lambda_* = \sum_{k=1}^{\infty}f(\x_k)-f(\x_{k+1})
      \geq \sum_{k=1}^{\infty}\eta\alpha_k\|\p(\x_k)\|^2
      \geq \eta\alpha_{\min} \sum_{k=1}^{\infty}\|\p(\x_k)\|^2.
  \end{equation*}
  It yields that
  \begin{eqnarray}\label{eq4-02}
    \sum_k \|\p(\x_k)\|^2 \leq \frac{2M}{\eta\alpha_{\min}}< +\infty.
  \end{eqnarray}
  Thus, the limit (\ref{eq4-05}) holds.

  Let $\x_{\infty}$ be an accumulation point of $\{\x_k\}$. Then
  $\x_{\infty}$ belongs to the compact set $\SPHERE$ and $\|\p(\x_{\infty})\|=0$.
  According to Lemma \ref{Lm3-2}, $\x_{\infty}$ is a generalized eigenvector
  whose associated eigenvalue is $f(\x_{\infty})=\lambda_*$.
\ep

\subsection{Further results based on the Kurdyka-{\L}ojasiewicz property}

In this subsection, we will prove that the iterates $\{\x_k\}$ generated by ACSA converge
without an assumption of the second-order sufficient condition.
The key tool of our analysis is the Kurdyka-{\L}ojasiewicz property.
This property was first discovered by S. {\L}ojasiewicz \cite{L63} in 1963 for real-analytic functions.
Bolte et al. \cite{BDL07} extended this property to nonsmooth subanalytic functions.
Whereafter, the Kurdyka-{\L}ojasiewicz property was widely applied to analyze
regularized algorithms for nonconvex optimization \cite{AB09,ABRS10}.
Significantly, it seems to be new to use the Kurdyka-{\L}ojasiewicz property to analyze
an inexact line search algorithm, e.g., ACSA proposed in Section 3.


We now write down the Kurdyka-{\L}ojasiewicz property \cite[Theorem 3.1]{BDL07} for completeness.

\begin{Theorem}[Kurdyka-{\L}ojasiewicz (KL) property]\label{Th4-05}
  Suppose that $\x_*$ is a critical point of $f(\x)$.
  Then there is a neighborhood $\mathscr{U}$ of $\x_*$,
  an exponent $\theta\in[0,1)$, and a constant $C_1$ such that
  for all $\x\in\mathscr{U}$, the following inequality holds
  \begin{equation}\label{eq4-09}
    \frac{|f(\x)-f(\x_*)|^{\theta}}{\|\p(\x)\|} \leq C_1.
  \end{equation}
  Here, we define $0^0\equiv1$.
\end{Theorem}

%

\begin{Lemma}\label{Lm4-06}
  Suppose that $\x_*$ is one of the accumulation points of $\{\x_k\}$.
  For the convenience of using the Kurdyka-{\L}ojasiewicz property,
  we assume that the initial iterate $\x_1$ satisfies
  $\x_1\in\mathscr{B}(\x_*,\rho)\equiv\{\x\in\Re^n~|~\|\x-\x_*\|<\rho\}\subseteq\mathscr{U}$ where
  \begin{equation*}
    \rho > \frac{2C_1}{\eta(1-\theta)}|f(\x_1)-f(\x_*)|^{1-\theta}+\|\x_1-\x_*\|.
  \end{equation*}
  Then, we have the following two assertions:
  \begin{equation}\label{eq4-10}
    \x_k\in\mathscr{B}(\x_*,\rho), \qquad \forall~k=1,2,\ldots,
  \end{equation}
  and
  \begin{equation}\label{eq4-11}
    \sum_k \|\x_{k+1}-\x_k\|\leq \frac{2C_1}{\eta(1-\theta)}|f(\x_1)-f(\x_*)|^{1-\theta}.
  \end{equation}
\end{Lemma}
\noindent
{\bf Proof}
  We prove (\ref{eq4-10}) by the induction.
  First, it is easy to see that $\x_1\in\mathscr{B}(\x_*,\rho)$.
  Next, we assume that there is an integer $K$ such that
  \begin{equation*}
    \x_k\in\mathscr{B}(\x_*,\rho), \qquad \forall~1 \leq k \leq K.
  \end{equation*}
  Hence, the KL property (\ref{eq4-09}) holds in these iterates.
  Finally, we now prove that $\x_{K+1}\in\mathscr{B}(\x_*,\rho)$.

  For the convenience of presentation, we define a scalar function
  \begin{equation*}
    \varphi(s)\equiv \frac{C_1}{1-\theta}|s-f(\x_*)|^{1-\theta}.
  \end{equation*}
  Obviously, $\varphi(s)$ is a concave function and its derivative is
  $\varphi'(s)=\frac{C_1}{|s-f(\x_*)|^{\theta}}$ if $s>f(\x_*)$.
  Then, for any $1\leq k \leq K$, we have
  \begin{eqnarray*}
    \varphi(f(\x_k))-\varphi(f(\x_{k+1}))
      &\geq& \varphi'(f(\x_k))(f(\x_k)-f(\x_{k+1})) \\
      & = & \frac{C_1}{|f(\x_k)-f(\x_*)|^{\theta}}(f(\x_k)-f(\x_{k+1})) \\
    \text{[by KL property]\qquad}  &\geq& \frac{1}{\|\p(\x_k)\|}(f(\x_k)-f(\x_{k+1})) \\
    \text{[since (\ref{suf-dec})]\qquad}  &\geq& \frac{1}{\|\p(\x_k)\|}\eta\alpha_k\|\p(\x_k)\|^2 \\
      &\geq& \frac{\eta\alpha_k\|\p(\x_k)\|}{\sqrt{1+\alpha_k^2\|\p(\x_k)\|^2}} \\
    \text{[because of (\ref{eq4-07})]\qquad}  &\geq& \frac{\eta}{2}\|\x_{k+1}-\x_k\|.
  \end{eqnarray*}
  It yields that
  \begin{eqnarray}\label{eq4-12}
    \sum_{k=1}^K \|\x_{k+1}-\x_k\|
      &\leq& \frac{2}{\eta}\sum_{k=1}^K \varphi(f(\x_k))-\varphi(f(\x_{k+1})) \nonumber\\
      &=& \frac{2}{\eta}(\varphi(f(\x_1))-\varphi(f(\x_{K+1}))) \nonumber\\
      &\leq& \frac{2}{\eta}\varphi(f(\x_1)).
  \end{eqnarray}
  So, we get
  \begin{eqnarray*}
    \|\x_{K+1}-\x_*\| &\leq& \sum_{k=1}^K\|\x_{k+1}-\x_k\|+\|\x_1-\x_*\| \\
      &\leq& \frac{2}{\eta}\varphi(f(\x_1))+\|\x_1-\x_*\| \\
      & < & \rho.
  \end{eqnarray*}
  Thus, $\x_{K+1}\in\mathscr{B}(\x_*,\rho)$ and (\ref{eq4-10}) holds.

  Moreover, let $K\to\infty$ in (\ref{eq4-12}). We obtain (\ref{eq4-11}).
\ep

\begin{Theorem}
  Suppose that the infinite sequence of iterates $\{\x_k\}$ is generated by ACSA.
  Then, the total sequence $\{\x_k\}$ has a finite length, i.e.,
  \begin{equation*}
    \sum_k \|\x_{k+1}-\x_k\| < +\infty,
  \end{equation*}
  and hence the total sequence $\{\x_k\}$ converges to a unique critical point.
\end{Theorem}
\noindent
{\bf Proof}
  Since the domain of $f(\x)$ is compact,
  the infinite sequence $\{\x_k\}$ generated by ACSA must have an accumulation point $\x_*$.
  According to Theorem \ref{Th4-04}, $\x_*$ is a critical point.
  Hence, there exists an index $k_0$, which could be viewed as an initial iteration when we use
  Lemma \ref{Lm4-06}, such that $\x_{k_0}\in\mathscr{B}(\x_*,\rho)$.
  From Lemma \ref{Lm4-06}, we have $\sum_{k=k_0}^{\infty} \|\x_{k+1}-\x_k\| < +\infty$.
  Therefore, the total sequence $\{\x_k\}$ has a finite length and converges to a unique critical point.
\ep

\begin{Lemma}
  There exists a positive constant $C_2$ such that
  \begin{equation}\label{eq4-13}
    \|\x_{k+1}-\x_k\| \geq C_2\|\p(\x_k)\|.
  \end{equation}
\end{Lemma}
\noindent
{\bf Proof}
  Since $\alpha_{\max}\geq\alpha_k\geq\alpha_{\min}>0$ and (\ref{eq4-07}), we have
  \begin{equation*}
    \|\x_{k+1}-\x_k\|=\frac{2\alpha_k\|\p(\x_k)\|}{\sqrt{1+\alpha_k^2\|\p(\x_k)\|^2}}
      \geq \frac{2\alpha_{\min}}{1+\alpha_{\max}M}\|\p(\x_k)\|.
  \end{equation*}
  Let $C_2\equiv\frac{2\alpha_{\min}}{1+\alpha_{\max}M}$. We get this lemma.
\ep

\begin{Theorem}\label{Th4-06}
  Suppose that $\x_*$ is the critical point of the infinite sequence of iterates $\{\x_k\}$
  generated by ACSA. Then, we have the following estimations.
  \begin{itemize}
    \item If $\theta\in(0,\tfrac{1}{2}]$, there exists a $\gamma>0$ and $\varrho\in(0,1)$ such that
      \begin{equation*}
        \|\x_k-\x_*\| \leq \gamma\varrho^k.
      \end{equation*}
    \item If $\theta\in(\tfrac{1}{2},1)$, there exists a $\gamma>0$ such that
      \begin{equation*}
        \|\x_k-\x_*\| \leq \gamma k^{-\frac{1-\theta}{2\theta-1}}.
      \end{equation*}
  \end{itemize}
\end{Theorem}
\noindent
{\bf Proof}
  Without loss of generality, we assume that $\x_1\in\mathscr{B}(\x_*,\rho)$.
  For convenience of following analysis, we define
  \begin{equation*}
    \Delta_k \equiv \sum_{i=k}^{\infty} \|\x_i-\x_{i+1}\| \geq \|\x_k-\x_*\|.
  \end{equation*}
  Then, we have
  \begin{eqnarray}
    \Delta_k &=& \sum_{i=k}^{\infty}\|\x_i-\x_{i+1}\| \nonumber\\
    \text{ [since (\ref{eq4-11})] }  &\leq& \frac{2C_1}{\eta(1-\theta)}|f(\x_k)-f(\x_*)|^{1-\theta} \nonumber\\
    &=& \frac{2C_1}{\eta(1-\theta)}\left(|f(\x_k)-f(\x_*)|^{\theta}\right)^{\frac{1-\theta}{\theta}} \nonumber\\
    \text{ [KL property] } &\leq&
        \frac{2C_1}{\eta(1-\theta)}\left(C_1\|\p(\x_k)\|\right)^{\frac{1-\theta}{\theta}} \nonumber\\
    \text{ [for (\ref{eq4-13})] } &\leq&
        \frac{2C_1}{\eta(1-\theta)}\left(C_1C_2^{-1}\|\x_k-\x_{k+1}\|\right)^{\frac{1-\theta}{\theta}} \nonumber\\
    &=& \frac{2C_1^{\frac{1}{\theta}}C_2^{-\frac{1-\theta}{\theta}}}{\eta(1-\theta)}
        \left(\Delta_k-\Delta_{k+1}\right)^{\frac{1-\theta}{\theta}} \nonumber\\
    &\equiv& C_3\left(\Delta_k-\Delta_{k+1}\right)^{\frac{1-\theta}{\theta}},  \label{eq4-14}
  \end{eqnarray}
  where $C_3$ is a positive constant.

  If $\theta\in(0,\tfrac{1}{2})$, we have $\tfrac{1-\theta}{\theta}\geq 1$.
  When the iteration $k$ is large enough, the inequality (\ref{eq4-14}) implies that
  \begin{equation*}
    \Delta_k \leq C_3(\Delta_k-\Delta_{k+1}).
  \end{equation*}
  That is
  \begin{equation*}
    \Delta_{k+1} \leq \frac{C_3-1}{C_3}\Delta_k.
  \end{equation*}
  Hence, recalling $\|\x_k-\x_*\|\leq\Delta_k$,
  we obtain the estimation if we take $\varrho\equiv\frac{C_3-1}{C_3}$.

  Otherwise, we consider the case $\theta\in(\tfrac{1}{2},1)$.
  Let $h(s)=s^{-\frac{\theta}{1-\theta}}$. Obviously, $h(s)$ is monotonously decreasing.
  Then, the inequality (\ref{eq4-14}) could be rewritten as
  \begin{eqnarray*}
    C_3^{-\frac{\theta}{1-\theta}} &\leq& h(\Delta_k)(\Delta_k-\Delta_{k+1}) \\
      &=& \int_{\Delta_{k+1}}^{\Delta_k} h(\Delta)~\mathrm{d}s \\
      &\leq& \int_{\Delta_{k+1}}^{\Delta_k} h(s)~\mathrm{d}s \\
      &=& -\frac{1-\theta}{2\theta-1}(\Delta_k^{-\frac{2\theta-1}{1-\theta}}
              -\Delta_{k+1}^{-\frac{2\theta-1}{1-\theta}}).
  \end{eqnarray*}
  Denote $\nu\equiv -\frac{1-\theta}{2\theta-1}<0$ since $\theta\in(\tfrac{1}{2},1)$.
  Then, we get
  \begin{equation*}
    \Delta_{k+1}^{\nu}-\Delta_k^{\nu} \geq \nu C_3^{-\frac{\theta}{1-\theta}} \equiv C_4>0.
  \end{equation*}
  It yields that for all $K>k$,
  \begin{equation*}
    \Delta_k \leq [\Delta_K^\nu+C_4(k-K)]^{\frac{1}{\nu}} \leq \gamma k^{\frac{1}{\nu}},
  \end{equation*}
  where the last inequality holds when the iteration $k$ is sufficiently large.
\ep

We remark that if the Hessian $H(\x_*)$ at the critical point $\x_*$ is positive definite,
the key parameter $\theta$ in the Kurdyka-{\L}ojasiewicz property is $\theta=\frac{1}{2}$.
Under Theorem \ref{Th4-06}, the sequence of iterates generated by ACSA has a linear convergence rate.
In this viewpoint, the Kurdyka-{\L}ojasiewicz property is weaker than the second order sufficient condition
of $\x_*$ being a minimizer.

\section{Numerical experiments}

To show the efficiency of the proposed ACSA algorithm, we perform some numerical experiments.
The parameters used in ACSA are
\begin{equation*}
    \eta=.001,\qquad \beta=.5,\qquad \alpha_{\max}=10000.
\end{equation*}
We terminate the algorithm if the objectives satisfy
$$\frac{|\lambda_{k+1}-\lambda_k|}{\max(1,|\lambda_k|)}< 10^{-12}\sqrt{n}$$
or the number of iterations exceeds $1000$.
The codes are written in MATLAB R2012a and run in a desktop computer with
Intel Core E8500 CPU at 3.17GHz and 4GB memory running Windows 7.

We will compare the following four algorithms in this section.
\begin{itemize}
  \item An adaptive shifted power method \cite{KM11,KM14} (Power M.) is implemented as
    {\tt eig\_sshopm} and {\tt eig\_geap} in Tensor Toolbox 2.6 for Z- and H-eigenvalues
    of even order symmetric tensors.
  \item An unconstrained optimization approach \cite{Han13} (Han's UOA) is solved by {\tt fminunc} in MATLAB
    with settings:  {\tt GradObj:on, LargeScale:off, TolX:1.e-10, TolFun:1.e-8,
    MaxIter:10000, Display:off.}
  \item For general symmetric tensors without considering a Hankel structure, we implement ACSA
    as ACSA-general.
  \item The ACSA algorithm (ACSA-Hankel) is proposed in Section 3 for Hankel tensors.
\end{itemize}

\subsection{Small Hankel tensors}

First, we examine some small tensors, whose Z- and H-eigenvalues could be computed exactly.


\begin{Example}[\cite{NW14}]\label{Exm-02}
  A Hankel tensor $\A$ whose entries are defined as
  \begin{equation*}
    a_{i_1i_2\cdots i_m}=\sin(i_1+i_2+\cdots+i_m), \qquad
      i_j=1,2,\ldots,n, ~j=1,2,\ldots,m.
  \end{equation*}
  Its generating vector is $\vv=(\sin(m), \sin(m+1), \ldots, \sin(mn))^{\top}$.

  If $m=4$ and $n=5$, there are five Z-eigenvalues which are listed as follows \cite{CDN,CHZ}
  \begin{equation*}
    \lambda_1=7.2595, ~~\lambda_2=4.6408, ~~\lambda_3=0.0000,
    ~~\lambda_4=-3.9204,~~\lambda_5=-8.8463.
  \end{equation*}
\end{Example}

\begin{table}[!bthp]
  \caption{Computed Z-eigenvalues of the Hankel tensor in Example \ref{Exm-02}.}\label{Tab-Two}
  \begin{center}
  \begin{tabular}{c|c|c|c|c}
    \hline
    Algorithms & Power M. & Han's UOA & ACSA-general & ACSA-Hankel \\
    \hline
     -8.846335 &  54\% & 58\% & 72\% &  72\% \\
     -3.920428 &  46\% & 42\% & 28\% &  28\% \\
    \hline
    CPU t. (sec)  & 23.09 & 9.34 & 8.39 & 0.67 \\
    \hline
  \end{tabular}
  \end{center}
\end{table}

We test four kinds of algorithms: power method, Han's UOA, ACSA-general and ACSA-Hankel.
For the purpose of obtaining the smallest Z-eigenvalue of the Hankel tensor,
we select $100$ random initial points on the unit sphere. The entries of each initial point
is first chosen to have a Gaussian distribution, then we normalize it to a unit vector.
The resulting Z-eigenvalues and CPU times are reported in Table \ref{Tab-Two}.
All of the four methods find the smallest Z-eigenvalue $-8.846335$.
But the occurrences for each method finding the smallest Z-eigenvalue are different.
We say that the ACSA algorithm proposed in Section 3 could find the extremal eigenvalues
with a higher probability.

Form the viewpoint of totally computational times, ACSA-general, and ACSA-Hankel are faster than
the power method and Han's UOA. When the Hankel structure of a fourth order
five dimensional symmetric tensor $\A$ is exploited, it is unexpected that
the new method is about 30 times faster than the power method.

\begin{Example}\label{Exm-03}
We study a parameterized fourth order four dimensional Hankel tensor $\mathcal{H}_{\epsilon}$
whose generating vector has the following form
\begin{equation*}
    \vv_{\epsilon}=(8-\epsilon,0,2,0,1,0,1,0,1,0,2,0,8-\epsilon)^\top.
\end{equation*}
If $\epsilon=0$, $\mathcal{H}_0$ is positive semidefinite but not positive definite \cite{CQW15}.
When the parameter $\epsilon$ is positive and trends to zero, the smallest Z- and H-eigenvalues
are negative and trends to zero.
In this example, we will illustrate this phenomenon by a numerical approach.
\end{Example}

\begin{figure}[tb!]
\begin{center}
  \includegraphics[scale=0.8]{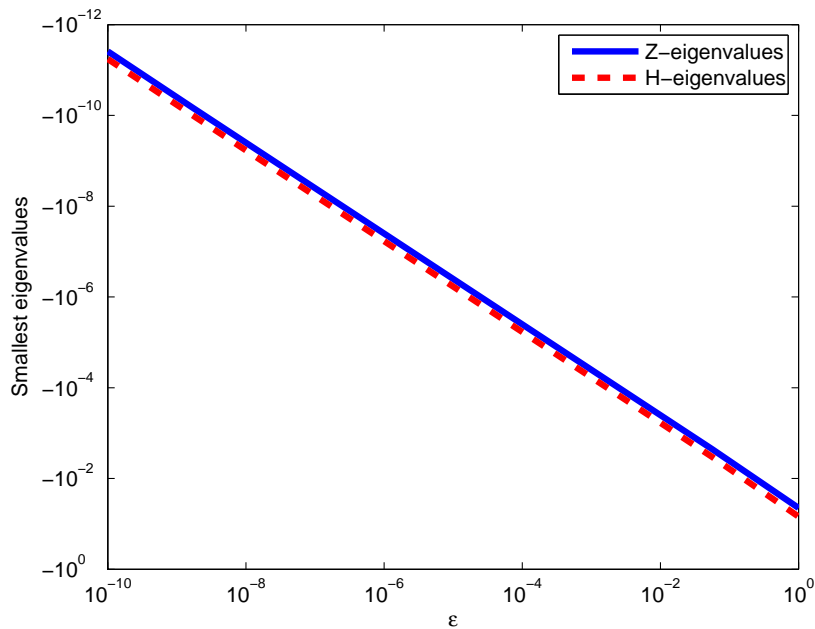}
\end{center}
\caption{The smallest Z- and H-eigenvalues of the parameterized fourth order four dimensional Hankel tensors.%
  }\label{Fig-Hankel}
\end{figure}

\begin{table}[!tbhp]
  \caption{CPU times (second) for computing Z- and H-eigenvalues of
    the parameterized Hankel tensors shown in Example \ref{Exm-03}.}\label{Tab-Six}
  \begin{center}
  \begin{tabular}{c|c|c|c|c}
    \hline
    Algorithms & Power M. & Han's UOA & ACSA-general & ACSA-Hankel \\
    \hline
     Z-eigenvalues   & 41.980 & 46.629 & 17.878 & 1.498 \\
     H-eigenvalues   & 29.562 & 45.833 & 16.973 & 1.544 \\
    \hline
     Total CPU times & 71.542 & 92.462 & 34.851 & 3.042 \\
    \hline
  \end{tabular}
  \end{center}
\end{table}

Again, we compare the power method, Han's UOA, ACSA-general, and ACSA-Hankel for
computing the smallest Z- and H-eigenvalues of the parameterized Hankel tensors in Example \ref{Exm-03}.
For the purpose of accuracy, we slightly modify the setting {\tt TolX:1.e-12, TolFun:1.e-12} for Han's UOA.
In each case, thirty random initial points on a unit sphere
are selected to obtain the smallest Z- or H-eigenvalues.
When the parameter $\epsilon$ decreases from $1$ to $10^{-10}$, the smallest Z- and H-eigenvalues returned
by these four algorithm are congruent. We show this results in Figure \ref{Fig-Hankel}.
When $\epsilon$ trends to zero, the smallest Z- and H-eigenvalues are negative and going to zero too.

The detailed CPU times for these four algorithms computing the smallest Z- and H-eigenvalues
of the parameterized fourth order four dimensional Hankel tensors are drawn in Table \ref{Tab-Six}.
Obviously, even without exploiting the Hankel structure, ACSA-general is two times faster than
the power method and Han's UOA. Furthermore, when the fast computational framework for
the products of a Hankel tensor time vectors is explored, ACSA-Hankel saves about $90\%$ CPU times.

\subsection{Large scale problems}

When the Hankel structure of higher order tensors is explored, we could compute
eigenvalues and associated eigenvectors of large scale Hankel tensors.

\begin{Example}\label{Exm-04}
  A Vandermonde tensor \cite{Qi15,Xu} is a special Hankel tensor.
  Let $$ \alpha=\frac{n}{n-1} \qquad\text{ and }\qquad \beta=\frac{1-n}{n}.$$
  Then, $\uu_1=(1,\alpha,\alpha^2,\ldots,\alpha^{n-1})^{\top}$
  and $\uu_2=(1,\beta,\beta^2,\ldots,\beta^{n-1})^{\top}$ are two Vandermonde vectors.
  The following $m$th order $n$ dimensional symmetric tensor
  \begin{equation*}
    \HH_V = \underbrace{ \uu_1\otimes \uu_1\otimes\cdots\otimes \uu_1}_{m\text{ times}}
        + \underbrace{ \uu_2\otimes \uu_2\otimes\cdots\otimes \uu_2}_{m\text{ times}}
  \end{equation*}
  is a Vandermonde tensor which satisfies the Hankel structure.
  Here $\otimes$ is the outer product. Obviously, the generating vector of $\HH_V$ is
  $\vv = (2,\alpha+\beta,\ldots,\alpha^{m(n-1)}+\beta^{m(n-1)})^{\top}$.
\end{Example}

\begin{Proposition}\label{Prop-9}
  Suppose the $m$th order $n$ dimensional Hankel tensor $\HH_V$ is defined as in Example \ref{Exm-04}.
  Then, when $n$ is even, the largest Z-eigenvalue of $\HH_V$ is $\|\uu_1\|^m$ and
  its associated eigenvector is $\tfrac{\uu_1}{\|\uu_1\|}$.
\end{Proposition}
\noindent
{\bf Proof}
  Since $\alpha\beta=-1$ and $n$ is even, $\uu_1$ and $\uu_2$ are orthogonal.
  We consider the optimization problem
  \begin{equation*}
  \begin{aligned}
    \max~&~ \HH_V\x^m = (\uu_1^{\top}\x)^m + (\uu_2^{\top}\x)^m, \\
    \mathrm{s.t.}~~&~\x^{\top}\x=1.
  \end{aligned}
  \end{equation*}
  Since $\|\uu_1\|>\|\uu_2\|$, when $\x=\tfrac{\uu_1}{\|\uu_1\|}$,
  the above optimization problem obtains its maximal value $\|\uu_1\|^m$.
  We write down its KKT condition, and it is easy to see that
  $(\|\uu_1\|^m, \tfrac{\uu_1}{\|\uu_1\|})$ is a Z-eigenpair of $\HH_V$.
\ep

\begin{table}[t!]
  \caption{The largest Z-eigenvalues of Vandermonde tensor in Example \ref{Exm-04}.%
    }\label{Tab-Three}
  \begin{center}
  \begin{tabular}{c|r|c|c|r}
    \hline
    $m$  & $n$ & largest Z-eigenvalues & Occurrences & CPU times (sec.) \\
    \hline
    4    &         10  & 9.487902e02   & 8 &     0.062 \\
    4    &        100  & 1.013475e05   & 8 &     0.140 \\
    4    &      1,000  & 1.019800e07   & 7 &     0.889 \\
    4    &     10,000  & 1.020431e09   & 8 &     9.048 \\
    4    &    100,000  & 1.020494e11   & 10 &  150.245 \\
    4    &  1,000,000  & 1.020500e13   & 5  & 2066.592 \\
    \hline
    6    &         10  & 2.922505e04    & 5 &     0.140 \\
    6    &        100  & 3.226409e07    & 5 &     0.234 \\
    6    &      1,000  & 3.256659e10    & 7 &     1.919 \\
    6    &     10,000  & 3.259683e13    & 7 &    17.753 \\
    6    &    100,000  & 3.259985e16    & 9 &   211.537 \\
    6    &  1,000,000  & 3.260016e19    & 4 &  3190.439 \\
    \hline
    8    &         10  & 9.002029e05    & 5 &     0.359 \\
    8    &        100  & 1.027131e10    & 5 &     0.437 \\
    8    &      1,000  & 1.039992e14    & 7 &     2.917 \\
    8    &     10,000  & 1.041279e18    & 7 &    30.561 \\
    8    &    100,000  & 1.041408e22    & 8 &  1058.248 \\
    \hline
  \end{tabular}
  \end{center}
\end{table}

Now, we employ the proposed ACSA algorithm which works with
the generating vector of a Hankel tensor to compute
the largest Z-eigenvalue of the Vandermonde tensor defined in Example \ref{Exm-04}.
We consider different orders $m=4,6,8$ and various dimension $n=10,\ldots,10^6$.
For each case, we choose ten random initial points, which has a Gaussian distribution
on a unit sphere. Table \ref{Tab-Three} shows the computed largest Z-eigenvalues
and the associated CPU times. For all case, the resulting largest Z-eigenvalue is agree with
Proposition \ref{Prop-9}. When the dimension of the tensor is one million,
the computational times for fourth order and sixth order Vandermonde tensors are
about $35$ and $55$ minutes respectively.

\begin{Example}\label{Exm-05}
  An $m$th order $n$ dimensional Hilbert tensor \cite{SQ14} is defined as
  \begin{equation*}
    \HH_{H}=\frac{1}{i_1+i_2+\cdots+i_m-m+1} \qquad
    i_j=1,2,\ldots,n, j=1,2,\ldots,m.
  \end{equation*}
  Its generating vector is
  $\vv=(1,\frac{1}{2},\frac{1}{3},\ldots,\frac{1}{m(n-1)+1})^{\top}$.
  When the order $m$ is even, the Hilbert tensors are positive definite.
  Its largest Z-eigenvalue and largest H-eigenvalues are bounded by
  $n^{\frac{m}{2}}\sin\frac{\pi}{n}$ and $n^{m-1}\sin\frac{\pi}{n}$ respectively.
\end{Example}

\begin{figure}[tb!]
\begin{center}
\begin{tabular}{cc}
  \includegraphics[scale=0.5]{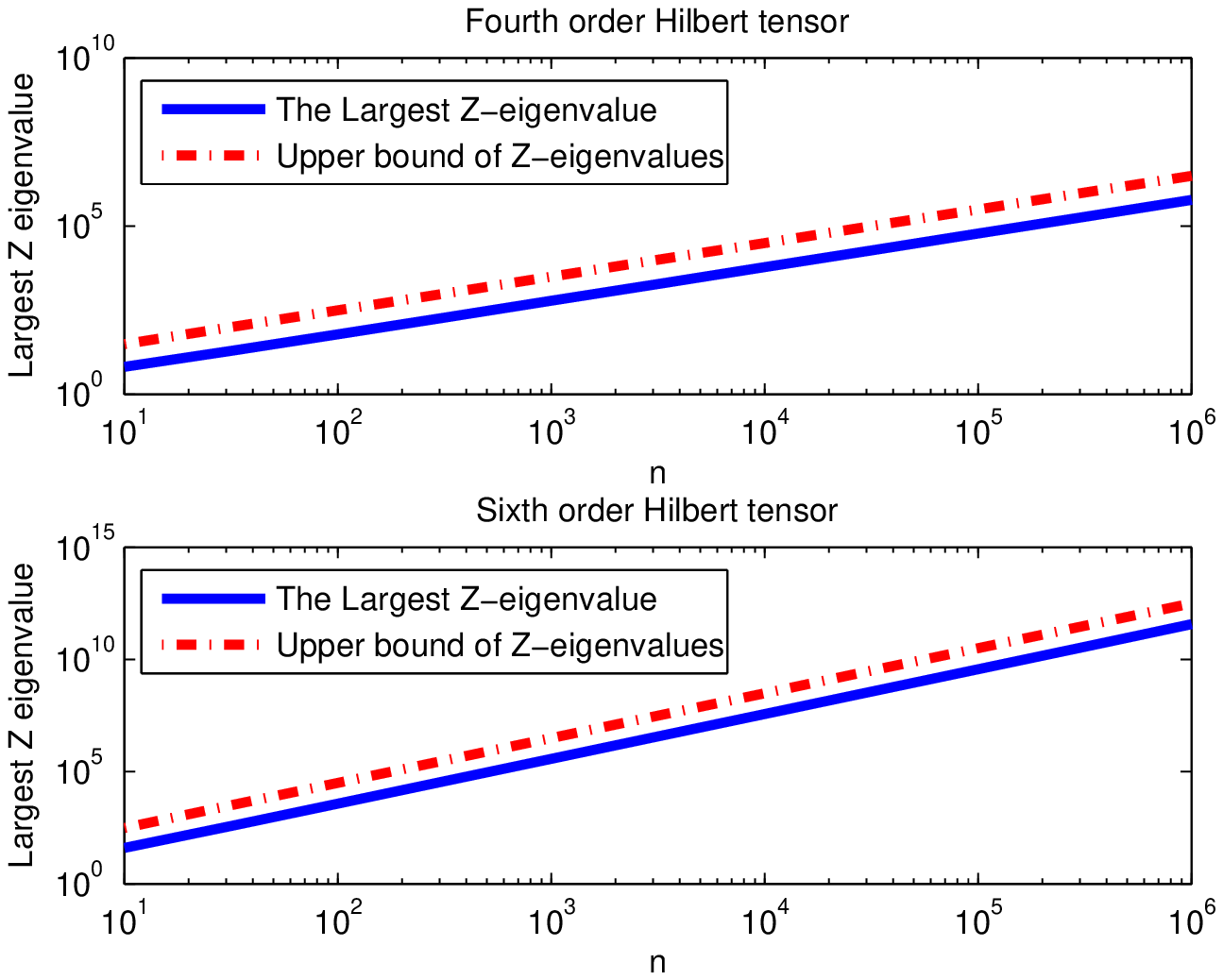} &
  \includegraphics[scale=0.5]{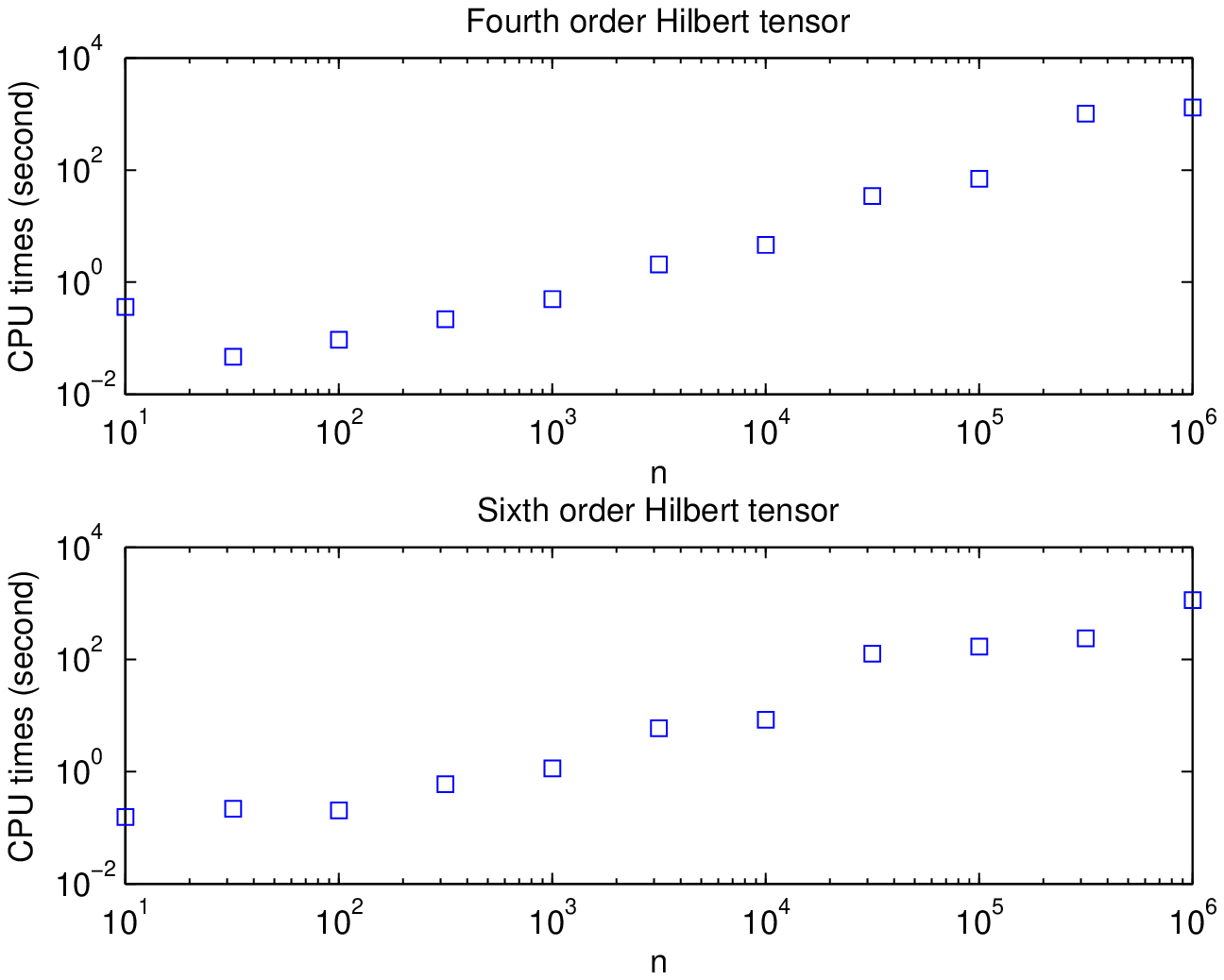} \\
\end{tabular}
\end{center}
\caption{The largest Z-eigenvalue and its upper bound for Hilbert tensors.%
  }\label{Fig-Hilbert-Zeigen}
\end{figure}

\begin{figure}[tb!]
\begin{center}
\begin{tabular}{cc}
  \includegraphics[scale=0.5]{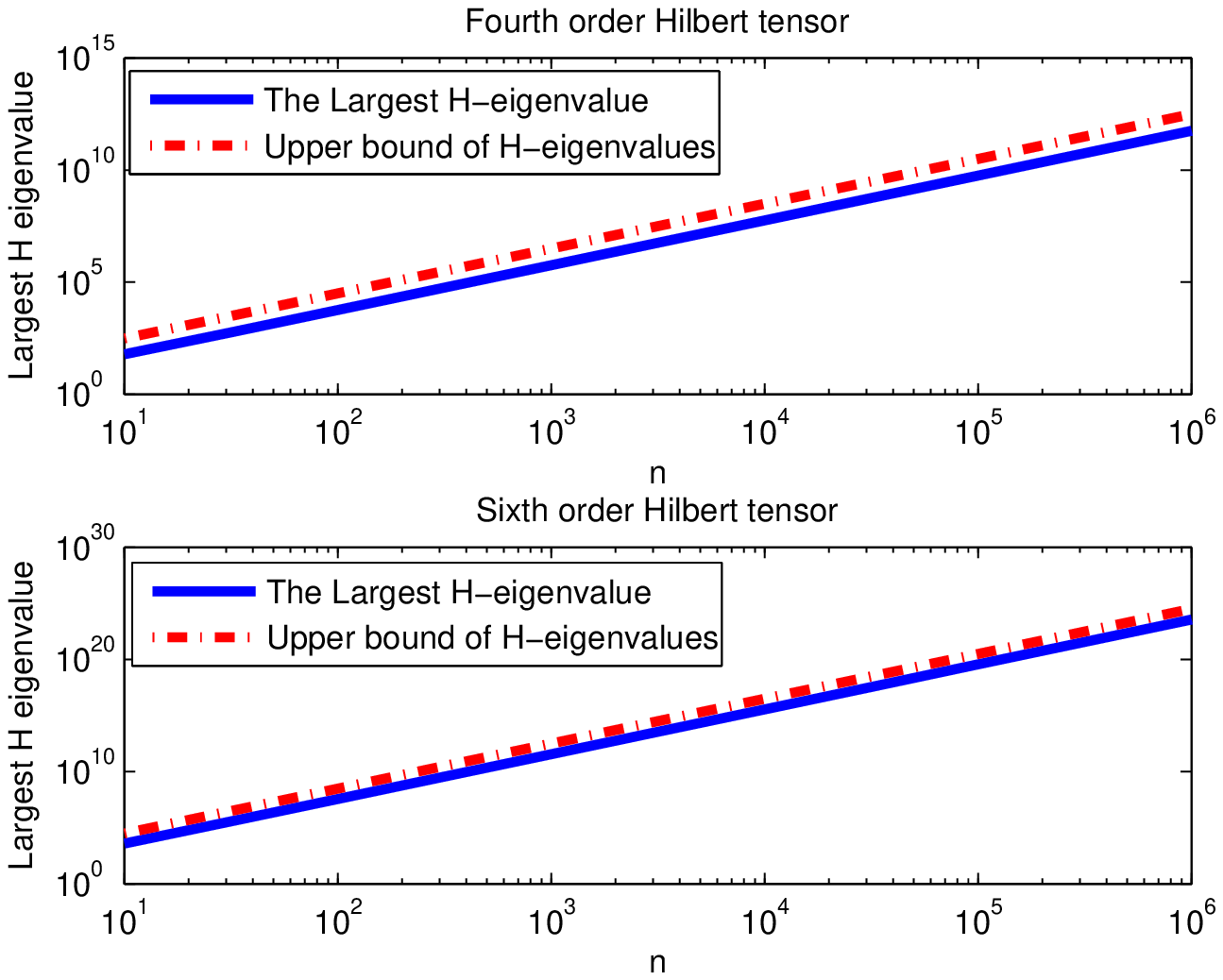} &
  \includegraphics[scale=0.5]{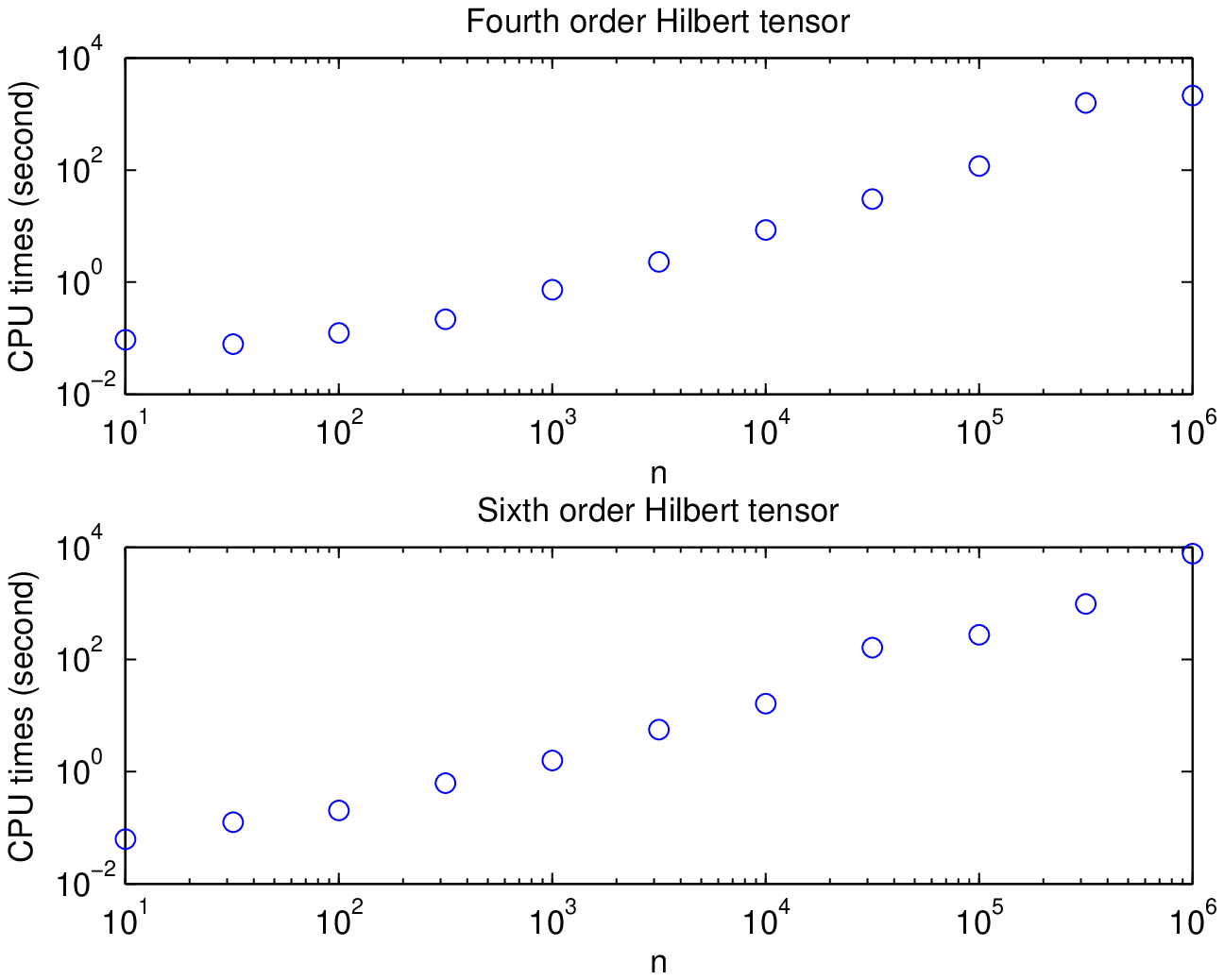} \\
\end{tabular}
\end{center}
\caption{The computed largest H-eigenvalue and its upper bound for Hilbert tensors.%
  }\label{Fig-Hilbert-Heigen}
\end{figure}

We illustrate by numerical experiments to show whether these bounds are tight?
First, for the dimension varying from ten to one million, we calculate the theoretical upper bounds
of the largest Z-eigenvalues of corresponding fourth order and sixth order Hilbert tensors.
Then, for each Hilbert tensor, we choose ten initial points and
employ the ACSA algorithm equipped with a fast computational framework for
products of a Hankel tensor and vectors to compute the largest Z-eigenvalues.
These results are shown in the left sub-figure of Figure \ref{Fig-Hilbert-Zeigen}.
The right sub-figure of Figure \ref{Fig-Hilbert-Zeigen} shows the corresponding CPU times
for ACSA-Hankel. We can see that the theoretical upper bounds for
the largest Z-eigenvalues of the Hilbert tensors are almost
tight up to a constant multiple.

Similar results for the largest H-eigenvalues and their theoretical upper bounds
of Hilbert tensors are illustrated in Figure \ref{Fig-Hilbert-Heigen}.

\section{Conclusion}

We proposed an inexact steepest descent method processing on a unit sphere for
generalized eigenvalues and associated eigenvectors of Hankel tensors.
Owing to the fast computation framework for the products of a Hankel tensor and vectors,
the new algorithm is fast and efficient as shown by some preliminary numerical experiments.
Since the Hankel structure is well-exploited, the new method could deal with some large scale
Hankel tensors, whose dimension is up to one million in a desktop computer.

\section*{Acknowledgment} The first author thanks Mr. Weiyang Ding and Dr. Ziyan Luo
  for the discussion on numerical experiments.


\begin{thebibliography}{99}
\bibitem{AB09} H. Attouch and J. Bolte, ``On the convergence of the proximal algorithm for nonsmooth functions
  involving analytic features'', {\sl Math. Program., Ser. B \bf 116} (2009) 5-16.

\bibitem{ABRS10} H. Attouch, J. Bolte, P. Redont and A. Soubeyran, ``Proximal alternating minimization and
  projection methods for nonconvex problems: an approach based on the Kurdyka-{\L}ojasiewicz inequality'',
  {\sl Math. Oper. Res. \bf 35} (2010) 438-457.

\bibitem{BK07} B. Bader and T. Kolda, ``Efficient MATLAB computations with sparse and factored tensors'',
  {\sl SIAM J. Sci. Comput. \bf 30} (2007) 205-231.

\bibitem{BB88} J. Barzilai and J.M. Borwein, ``Two-point step size gradient methods'',
  {\sl IMA J. Numer. Anal., \bf 8}, (1988) 141-148.

\bibitem{BDL07} J. Bolte, A. Daniilidis and A. Lewis, ``The {\L}ojasiewicz inequality for nonsmooth subanalytic
  functions with applications to subgradient dynamical systems'', {\sl SIAM J. Optim. \bf 17}, (2006) 1205-1223.

\bibitem{BDA07} R. Boyer, L. De Lathauwer and K. Abed-Meraim, ``Higher order tensor-based method for
  Delayed exponential fitting'', {\sl IEEE T. Signal Proces. \bf 55} (2007) 2795-2809.

\bibitem{CPZ09} K.C. Chang, K. Pearson and T. Zhang, ``On eigenvalue problems of real symmetric tensors'',
  {\sl J. Math. Anal. Appl. \bf 350} (2009) 416-422.

\bibitem{CHZ} L. Chen, L. Han and L. Zhou, ``Computing tensor eigenvalues via homotopy methods'',
  (2015) ``http://arxiv.org/pdf/1501.04201v3.pdf''.

\bibitem{CDHS13} Y. Chen, Y. Dai, D. Han, and W. Sun, ``Positive semidefinite generalized diffusion tensor
  imaging via quadratic semidefinite programming'', {\sl SIAM J. Imaging Sci., \bf 6} (2013) 1531-1552.

\bibitem{CQW15} Y. Chen, L. Qi and Q. Wang, ``Positive semi-definiteness and sum-of-squares property of
  fourth order four dimensional Hankel tensors'', (2015) ``http://arxiv.org/pdf/1502.04566v8.pdf''.

\bibitem{CV14} J.H. Choi and S.V.N. Vishwanathan, ``DFacTo: distributed factorization of tensors'',
  (2014) ``http://arxiv.org/pdf/1406.4519v1.pdf''.

\bibitem{CiP09} A. Cichocki and A.-H. Phan, ``Fast Local algorithms for large scale nonnegative matrix
  and tensor factorizations'', {\sl IEICE T. Fund. Electr. \bf E92-A} (2009) 708-721.

\bibitem{CD12} J. Cooper and A. Dutle, ``Spectra of uniform hypergraphs'',
  {\sl Linear algebra Appl. \bf 436} (2012) 3268-3292.

\bibitem{CDN} C. Cui, Y. Dai and J. Nie, ``All real eigenvalues of symmetric tensors'',
  {\sl SIAM J. Matrix Anal. Appl. \bf 35} (2014) 1582-1601.

\bibitem{Dai14} Y. Dai, ``A positive BB-like stepsize and an extension for symmetric linear systems'',
  in {\sl Workshop on Optimization for Modern Computation}, Beijing, China, (2014),
  ``http://bicmr.pku.edu.cn/conference/opt-2014/slides/Yuhong-Dai.pdf''.

\bibitem{DK14} A.L.F. de Almeida and A.Y. Kibangou, ``Distributed large-scale tensor decomposition'',
  in {\sl IEEE International Conference on Acoustics, Speech and Siganl Processing (ICASSP)} (2014) 26-30.

\bibitem{DDV00} L. De Lathauwer, B. De Moor, and J. Vandewalle,
  ``On the best rank-$1$ and rank-$(R_1,R_2,\ldots,R_N)$ approximation of higher-order tensors'',
  {\sl SIAM J. Matrix Anal. Appl. \bf 21} (2000) 1324-1342.

\bibitem{DQW} W. Ding, L. Qi and Y. Wei,
  ``Fast Hankel tensor-vector product and its application to exponential data fitting'',
  {\sl Numer. Linear Algebr. Appl.}, (2015), DOI: 10.1002/nla.1970.

\bibitem{DW15} W. Ding and Y. Wei, ``Generalized tensor eigenvalue problems'',
  {\sl SIAM J. Matrix Anal. Appl.} (2015), To appear.

\bibitem{FNO87} S. Friedland, J. Nocedal and M.L. Overton, ``The formulation and analysis of numerical
  methods for inverse eigenvalue problems'', {\sl SIAM J. Numer. Anal. \bf 24} (1987) 634-667.

\bibitem{GWY09} D. Goldfarb, Z. Wen and W. Yin, ``A curvilinear search method for the $p$-harmonic flow
  on spheres'', {\sl SIAM J. Imaging Sci., \bf 2} (2009) 84-109.

\bibitem{GV13} G.H. Golub and C.F. Van Loan, {\sl Matrix Computations, 4th Edition},
  The Johns Hopkins University Press (2013), ISBN 978-1-4214-0794-4.

\bibitem{Han13} L. Han, ``An unconstrained optimization approach for finding real eigenvalues of
  even order symmetric tensors'', {\sl Numer. Algebr. Control Optim. \bf 3} (2013) 583-599.

\bibitem{HCD} C. Hao, C. Cui and Y. Dai, ``A sequential subspace projection method for extreme Z-eigenvalues
  of supersymmetric tensors'', {\sl Numer. Linear Algebr. Appl. \bf 22} (2015) 283-298.

\bibitem{HCD15a} C. Hao, C. Cui and Y. Dai,  ``A feasible trust-region method for calculating extreme
Z-eigenvalues of symmetric tensors'',  {\sl Pacific J. Optim.}, in press (2015).

\bibitem{HL13} C.J. Hillar and L.-H. Lim, ``Most tensor problems are NP-hard'', {\sl J. ACM \bf 60} (2013)
  article 45:1-39.

\bibitem{HHQ} S. Hu, Z. Huang and L. Qi, ``Finding the extreme Z-eigenvalues of tensors via
  a sequential SDPs method'', {\sl Numer. Linear Algebr. Appl. \bf 20} (2013) 972-984.

\bibitem{KPHF12} U. Kang, E. Papalexakis, A. Harpale, and C. Faloutsos, ``GigaTensor: scaling tensor
  analysis up by 100 times - algorithms and discoveries'', in {\sl Proceedings of the 18th ACM SIGKDD
  International Conference on Knowledge Discovery and Data Mining}, (2012) 316-324.

\bibitem{KR02} E. Kofidis and P.A. Regalia,
  ``On the best rank-$1$ approximation of higher-order supersymmetric tensors'',
  {\sl SIAM J. Matrix Anal. Appl. \bf 23} (2002) 863-884.

\bibitem{KM11} T.G. Kolda and J.R. Mayo, ``Shifted power method for computing tensor eigenpairs'',
  {\sl SIAM J. Matrix Anal. Appl. \bf 32} (2011) 1095-1124.

\bibitem{KM14} T.G. Kolda and J.R. Mayo,
  ``An adaptive shifted power method for computing generalized tensor eigenpairs'',
  {\sl SIAM J. Matrix Anal. Appl. \bf 35} (2014) 1563-1581.

\bibitem{Lim} L.-H. Lim, ``Singular values and eigenvalues of tensors: a variational approach'', in
  {\sl Proceedings of the IEEE International Workshop on Computational Advances in Multi-Sensor Adaptive Processing
  (CAMSAP'05), \bf 1} (2005) 129-132.

\bibitem{L63} S. {\L}ojasiewicz, ``Une propri\'{e}t\'{e} topologique des sous-ensembles analytiques r\'{e}els'',
  {\sl Les \'{E}quations aux D\'{e}riv\'{e}es Partielles}, \'{E}ditions du centre National de la Recherche
  Scientifique, Paris, 87-89, (1963).

\bibitem{LT03} J.-G. Luque and J.-Y. Thibon, ``Hankel hyperdeterminants and Selberg integrals'',
  {\sl J. Phys. A. \bf 36} (2003) 5267-5292.

\bibitem{ML13} J. McAuley and J. Leskovec, ``Hidden factors and hidden topics: understanding rating
  dimensions with review text'', in {\sl Proceeding of the 7th ACM Conference on Recommender Systems},
  (2013) 165-172.

\bibitem{NQB14} G. Ni, L. Qi and M. Bai, ``Geometric measure of entanglement and U-eigenvalues of tensors'',
  {\sl SIAM J. Matrix Anal. Appl. \bf 35} (2014) 73-87.

\bibitem{NQW08} Q. Ni, L. Qi and F. Wang, ``An eigenvalue method for testing positive definiteness of
  a multivariate form'', {\sl IEEE T. Automat. Contr. \bf 53} (2008) 1096-1107.

\bibitem{NW14} J. Nie and L. Wang, ``Semidefinite relaxations for best rank-1 tensor approximations'',
  {\sl SIAM J. Matrix Anal. Appl. \bf 35} (2014) 1155-1179.

\bibitem{OS11} V. Oropeza and M. Sacchi, ``Simultaneous seismic data denoising and reconstruction via
  multichannel singular spectrum analysis'', {\sl Geophysics \bf 76} (2011) V25-V32.

\bibitem{PDV05} J.M. Papy, L. De Lathauwer and S. Van Huffel, ``Exponential data fitting using
  multilinear algebra: the single-channel and multi-channel case'', {\sl Numer. Linear Algebr. Appl. \bf 12}
  (2005) 809-826.

\bibitem{PDV09} J.M. Papy, L. De Lathauwer and S. Van Huffel, ``Exponential data fitting using
  multilinear algebra: the decimative case'', {\sl J. Chemometr. \bf 23} (2009) 341-351.

\bibitem{Qi05} L. Qi, ``Eigenvalues of a real supersymmetric tensor'',
  {\sl J. Symb. Comput. \bf 40} (2005) 1302-1324.

\bibitem{Qi15} L. Qi, ``Hankel tensors: Associated Hankel matrices and Vandermonde decomposition'',
  {\sl Commun. Math. Sci. \bf 13} (2015) 113-125.

\bibitem{QWW09} L. Qi, F. Wang and Y. Wang, ``Z-eigenvalue methods for a global polynomial optimization problem'',
  {\sl Math. Program., Ser. A \bf 118} (2009) 301-316.

\bibitem{QYX13} L. Qi, G. Yu and Y. Xu, ``Nonnegative diffusion orientation distribution function'',
  {\sl J. Math. Imaging. Vis. \bf 45} (2013) 103-113.

\bibitem{SLVK14} M.D. Schatz, T.-M. Low, R.A. Van De Geijn, and T.G. Kolda,
  ``Exploiting symmetry in tensors for high performance'',
  {\sl SIAM J. Sci. Comput. \bf 36} (2014) C453-C479.

\bibitem{SS08} T. Schultz and H.-P. Seidel, ``Estimating crossing fibers: a tensor decomposition approach'',
  {\sl IEEE T. Vis. Comput. Gr. \bf 14} (2008) 1635-1642.

\bibitem{S14} R.S. Smith, ``Frequency domain subspace identification using nuclear norm minimization and
  Hankel matrix realizations'', {\sl IEEE T. Automat. Contr. \bf 59} (2014) 2886-2896.

\bibitem{SQ14} Y. Song and L. Qi, ``Infinite and finite dimensional Hilbert tensors'',
  {\sl Linear Algebra. Appl. \bf 451} (2014) 1-14.

\bibitem{TBM13} S. Trickett, L. Burroughs and A. Milton, ``Interpolating using Hankel tensor completion'',
  in {\sl SEG Annual Meeting} (2013) 3634-3638.

\bibitem{VCDV} S. Van Huffel, H. Chen, C. Decanniere, and P. Van Hecke, ``Algorithm for time-domain
  NMR data fitting based on total least squares'', {\sl J. Magn. Reson., Series A \bf 110}
  (1994) 228-237.

\bibitem{WY} Z. Wen and W. Yin, ``A feasible method for optimization with orthogonality constraints'',
  {\sl Math. Program., Ser A \bf 142} (2013) 397-434.

\bibitem{Xu} C. Xu, ``Hankel tensors, Vandermonde tensors and their positivities'',
  {\sl Linear Algebra Appl.}, in press (2015).


\end{thebibliography}
\end{document}